\documentclass[11pt]{article}
\usepackage[T1]{fontenc}
\usepackage{lmodern,amsmath,amssymb,amsthm,mathtools}
\usepackage[margin=1in]{geometry}
\usepackage{enumitem,microtype,needspace}
\usepackage[colorlinks=true,linkcolor=blue,citecolor=blue,urlcolor=blue]{hyperref}
\usepackage[T1]{fontenc}
\usepackage{lmodern,amsmath,amssymb,amsthm,mathtools}
\usepackage[margin=1in]{geometry}
\usepackage{microtype,booktabs}

\newtheorem{theorem}{Theorem}[section]
\newtheorem{proposition}[theorem]{Proposition}
\newtheorem{lemma}[theorem]{Lemma}

\theoremstyle{definition}
\newtheorem{definition}[theorem]{Definition}

\numberwithin{equation}{section}
\newtheorem{question}{Question}
\newcommand{\normal}[1]{\langle\!\langle #1\rangle\!\rangle}
\newcommand{\Prob}{\mathbb P}
\newcommand{\Conj}{\operatorname{Conj}}
\newcommand{\rk}{\operatorname{rk}}
\newcommand{\eps}{\varepsilon}
\setlist{itemsep=3pt,topsep=5pt}
\title{Random Burnside groups}
\author{Olga Kharlampovich, Alexei Miasnikov}
\date{21 October 2026}
\begin{document}
\maketitle
\begin{abstract}
Fix integers $m\ge2$ and $k\ge1$. For every sufficiently large fixed odd exponent $M$, choose $k$ independent elements uniformly from the sphere of radius $\ell$ in the free Burnside group $B_m(M)$. We prove that the quotient by their normal closure is infinite with probability at least $1-e^{-c\ell}$. More precisely, the quotient map has no nontrivial kernel element of length at most $\kappa\ell$, for constants $c,\kappa>0$.  We also consider the analogue of the Gromov positive density model for the Burnside variety. \end{abstract}

\section{Statement and framework}
Random quotients provide a natural way to study which properties of a
group survive the addition of relations. For free Burnside groups,
the first question is whether the resulting quotient remains infinite.
Coulon and Gruber \cite{CoulonGruber} developed small cancellation
theory in Burnside varieties. Under suitable small cancellation
conditions and restrictions on powers occurring in the relators,
their results produce infinite groups of sufficiently large odd
exponent. Gruber and Mackay \cite{GruberMackay} subsequently studied
random triangular Burnside groups. They proved that, for each density
\[
  0<d<\frac{11-\sqrt{41}}{12},
\]
and every sufficiently large fixed exponent, these groups are infinite
with probability tending to one as the number of generators increases. Our model differs from the triangular model of Gruber and Mackay in both the limiting process and the sampling measure. In their model, the exponent is fixed, the number (m) of generators tends to infinity, and one chooses $(\lfloor(2m-1)^{3d}\rfloor)$ independent, uniformly distributed cyclically reduced words of length three in the free group, where (d) is a fixed density parameter. These words are then imposed as relations in the corresponding free Burnside group. In our model, the rank (m), a sufficiently large odd exponent (M), and the number (k) of additional relations are all fixed.

We study random quotients of a fixed free Burnside group
$B_m(M)$. We fix the rank $m$, a sufficiently large odd exponent $M$,
and the number $k$ of additional relations. These relations are chosen
independently and uniformly from the  sphere of radius $\ell$
in $B_m(M)$ (intrinsic sphere), and $\ell$ tends to infinity. We prove that the quotient
is infinite with exponentially high probability. More precisely, on
the same event, the quotient map has no nontrivial kernel element of
length at most $\kappa\ell$, for a constant $\kappa>0$. The proof
combines conjugacy-class counting with a graded diagram construction
in which the random cyclic cores are inserted at an intermediate
rank and all subsequent periods are reselected.

Let $F_m=F(x_1,\ldots,x_m)$ and
\[
 B=B_m(M)=F_m/\normal{W^M:W\in F_m}.
\]
All group lengths are taken with respect to the images of $x_1,\ldots,x_m$. Write $S_B(r)$ and $B_B(r)$ for the  sphere and closed ball of radius $r$, respectively. A word length $|W|$ is its literal length; $|g|_B$ is the intrinsic length of a group element. A \emph{cyclic core} of $g$ is a word of minimum length among all words representing conjugates of $g$ in $B$. Every cyclic shift of a cyclic core is geodesic in $B$.

\begin{theorem}\label{main}
For every $m\ge2$ and $k\ge1$ there is $M_0$ such that the following holds for every fixed odd $M\ge M_0$. There are $c,\kappa>0$ such that, for all sufficiently large $\ell$, independent uniform elements $g_1,\ldots,g_k\in S_B(\ell)$ satisfy
\[
 \Pr\!\left\{\ker(B\longrightarrow Q_\ell)\cap B_B(\kappa\ell)=\{1\}\right\}
 \ge1-e^{-c\ell},\qquad
 Q_\ell=B/\normal{g_1,\ldots,g_k}.
\]
In particular, $Q_\ell$ is infinite with probability at least $1-e^{-c\ell}$.
\end{theorem}

The sampling is uniform on group elements in a sphere, not on freely reduced words in $F_m$. The exponent is fixed before $\ell$ tends to infinity. The theorem includes sufficiently large odd composite exponents.

The proof has two parts. First, a finite list of compressed descriptions excludes all long overlaps through short connecting paths that will be needed in diagrams, uniformly over cyclic cores and their spellings. Second, a graded construction inserts these cores below their length and reselects all later periods. The diagram induction uses the excluded descriptions only after lowering each relevant mixed strip to an old rank. The final kernel estimate follows from the perimeter inequality for graded diagrams.

In the last section   we consider the analogue of the Gromov positive density model for the Burnside variety (Theorem \ref{gromovtheorem}). 

\section{Known facts and diagram conventions}\label{inputs}
We use the standard graded Burnside presentation of \cite[Sections 18--19]{Olshanskii}. Let $G_j$ be its quotient through integer rank $j$. At rank $j$, one selects representatives, up to conjugacy and inversion in $G_{j-1}$, of the words of length $j$ that are simple there, and imposes $A^M=1$ for each selected period $A$. A word is \emph{simple} at a stage if it is nontrivial and is not conjugate at that stage to a power of a shorter word or of a previously selected period. The direct limit of the $G_j$ is $B$. We call this presentation and its ranks \emph{old}.

\subsection{Facts about the old groups}
We use the following established results.
\begin{enumerate}[label=(\roman*),leftmargin=*]
\item The logarithmic growth rate $h_M$ of $B_m(M)$ tends to $\log(2m-1)$ as odd $M\to\infty$; apply \cite[Theorem 1.2]{Coulon} to $F_m$.
\item Short equalities and conjugacies in $B$ hold at old rank $O(r/M)$ when the words concerned have length at most $r$. We use the explicit stabilization bound in \cite[Lemma 5.1]{GGO}; in particular a bound $\lceil 3r/M\rceil$ suffices.
\item If $A$ is simple in an old stage, every $A$-periodic boundary section in a reduced diagram of that stage is $(1-\beta_o)$-quasigeodesic in the boundary comparison sense: a complementary path to a section $p$ has length at least $(1-\beta_o)|p|$. The parameter $\beta_o$ can be prescribed sufficiently small before choosing $M$ \cite[Lemma 4.9]{GGO}.
\item The centralizer in $B$ of a nonidentity element is cyclic of order $M$ \cite[Theorem 3.2]{Adian}. In particular $XY=YZ$ has at most $M$ solutions $Y$ when $X\ne1$ and $X,Z$ are fixed. Also, in every group of odd exponent $M$, $Y^2=Z$ has exactly one solution, namely $Z^{(M+1)/2}$.
\end{enumerate}
The quantitative conjugator estimate from \cite[Lemma 8.9]{GGO} is stated and used in Section~\ref{counting}. All choices of its small parameters precede the choice of the exponent.

\subsection{Graded maps}
We use disc and annular diagrams with positive cells assigned ranks. Edges labelled by the identity and the usual trivial cells are allowed and have length zero. Adding them to separate coincident arcs is called a \emph{zero-refinement}; it changes neither labels in the group nor the number of positive cells. A path is geodesic in a diagram if no path homotopic to it relative to its endpoints is shorter. The geometric arguments use only order and inequalities between positive ranks; we may therefore insert a half-integer rank between two consecutive integer ranks.

Contiguity subdiagrams are the recursively defined graded strips of \cite[Section 14]{Olshanskii}, also described in \cite[Section 2.3]{OS}. A zero-bond is the zero-width adjacency between positive edges. A positive bond uses a principal cell of smaller rank together with two previously defined strips covering at least an $\eps$-fraction of its perimeter. A strip between two cells, between a cell and a boundary section, or between two boundary sections is bounded by one such bond or two disjoint bonds. It has boundary $c p d q$, where $p,q$ lie on the indicated cells or sections and $c,d$ are the side arcs. Its degree from a source cell $\Pi$ is $|p|/|\partial\Pi|$. Boundary sections impose no additional rank inequality in this definition. Thus a strip remains a boundary-to-boundary strip after its endpoint cells are removed.

Fix a geometric perimeter parameter $n_0$, distinct from $M$. The following definitions and geometric estimates use $n_0$.
\Needspace{8\baselineskip}
\begin{definition}\label{amapdef}
A graded map satisfies \emph{Condition A} if:
\begin{enumerate}[label=A\arabic*.,leftmargin=*]
\item every positive cell $\Pi$ has $|\partial\Pi|\ge n_0\rk\Pi$;
\item every cell-contour interval of length at most $\max(\rk\Pi,2)$ is geodesic in the map;
\item in a contiguity strip of degree at least $\eps$ from a cell to a cell of rank $a$, the target arc has length less than $(1+\gamma)a$.
\end{enumerate}
Such a map is an \emph{A-map}. A boundary section $q$ is \emph{smooth of rank $a$} if its intervals of length at most $\max(a,2)$ are geodesic and the same target-arc bound holds for every cell-to-$q$ strip of degree at least $\eps$.
\end{definition}

We will use the following consequences of Condition A. They are statements about graded maps and do not require all positive-cell labels to be powers; see \cite[Lemmas 2.1, 2.7--2.8, 2.12, 2.14--2.16]{OS} and \cite[Sections 15--17]{Olshanskii}.
\begin{enumerate}[label=(G\arabic*),leftmargin=*]
\item\label{geom:sides} In a boundary-to-boundary contiguity strip in an A-map, if the two sections are smooth of ranks $a,b$, each side arc has length less than $\zeta n_0\min(a,b)$.
\item\label{geom:length} If $q$ is smooth, an alternative boundary path $v$ in a disc satisfies $(1-\beta)|q|\le|v|$. For an annulus the same inequality compares a smooth full contour with the other contour.
\item\label{geom:perimeter} If a disc A-map contains a positive cell $\Pi$, then $|\partial\Delta|>(1-\beta)|\partial\Pi|$.
\item\label{geom:rank} If a disc A-map has two smooth sides of ranks $a,b$ and each remaining side is shorter than $\zeta n_0\min(a,b)$, then its rank is less than $3\gamma^{-1}\zeta\min(a,b)$.
\item\label{geom:trim} Suppose the smooth sides have ranks $a,b$, put $s=\min(a,b)$, and let the two other sides have total length $C$. If both smooth sides are at least $T=\gamma^{-1}(C+s)$ long, a subdiagram can be cut out whose connectors have length less than $\alpha s$ and whose two smooth sides each lose less than $T$.
\item\label{geom:annulus} An annular A-map with no essential loop of length zero has a connecting path of length less than $\gamma$ times the sum of its contour lengths.
\end{enumerate}
We also use the interior-point cutting statement of \cite[Lemma 2.12]{OS}. Its parameters, and those in (G1)--(G6), can be chosen in the standard order $n_0^{-1}\ll\zeta\ll\eps\ll\delta_g\ll\gamma\ll\beta\ll\alpha\ll1$. All later comparisons keep these same constants.

\section{Counting cyclic cores and contacts}\label{counting}
Put $s(r)=|S_B(r)|$ and
\[
 h=h_M=\lim_{r\to\infty}\frac{\log s(r)}r,
 \qquad E(t)=\max_{0\le r\le t}(\log s(r)-hr),
 \qquad H=\log(2m-1),\quad D=\log(2m+2).
\]
Splitting a chosen geodesic gives $s(a+b)\le s(a)s(b)$, so subadditivity lemma gives $s(r)\ge e^{hr}$ and $E(t)=o(t)$. Here $s(0)=1$. Since $|B_B(r)|=\sum_{j=0}^r s(j)$, its logarithmic growth rate is the same $h$: for every $\eta>0$ one has $s(j)\le C_\eta e^{(h+\eta)j}$, while $|B_B(r)|\ge s(r)\ge e^{hr}$. Thus the quoted growth theorem also applies to the logarithmic growth rate of the spheres.  We take $M$ sufficiently large that $h>D/2$.

Fix a small coding parameter $\gamma_c$. For a nontrivial conjugacy class with minimum length $L$, the conjugator code of \cite[Lemma 8.9]{GGO} gives
\begin{equation}\label{fiber}
 |\Conj(U)\cap S_B(\ell)|
 \le L\exp\!\left(D\left[\frac{\ell-L}{2}+9\gamma_c\ell+1\right]\right).
\end{equation}
Indeed, fix one shortest word $U$ and fix the outer basepoint in a reduced conjugacy annulus. The code for the conjugator has length at most $(\ell-L)/2+9\gamma_c\ell$ over an alphabet of size $2m+2$. Its value and one of the $L$ cyclic shifts of $U$ determine the conjugate element. No count of all geodesic spellings of that element is involved.

\begin{lemma}[Long cyclic cores]\label{cores}
For $0<b<1$, if $9\gamma_cD<(1-b)(h-D/2)$, the probability that the minimum conjugacy length of a uniform sphere element is at most $b\ell$ is bounded by
\[
 \exp\!\left(-[(1-b)(h-D/2)-9\gamma_cD]\ell+o(\ell)\right).
\]
\end{lemma}
\begin{proof}
There are at most $s(L)$ conjugacy classes of minimum length $L$. Multiply by \eqref{fiber}, sum over $L\le b\ell$, and divide by $s(\ell)$. Polynomial factors and $E(\ell)$ contribute $o(\ell)$ to the exponent.
\end{proof}

\begin{lemma}\label{transfer}
Fix one of the finite decoding schemes below. Suppose an exceptional $r$-tuple of conjugacy classes, $r\in\{1,2\}$, has minimum lengths $L_1,\ldots,L_r\le\ell$ and can be reconstructed from \(a\) group elements, together with the auxiliary data specified below of total intrinsic length at most
\[
 \sum_i L_i-u.
\]
Suppose the auxiliary choices and decoding multiplicity are bounded by a polynomial in $\ell$ (with $M$ fixed). Then its probability for independent sphere inputs is at most
\begin{equation}\label{transferbound}
 \operatorname{poly}(\ell)
 \exp\!\left(-hu+9r\gamma_cD\ell+aE(r\ell)\right).
\end{equation}
The description may use any shortest representative and any cyclic basepoint. A recorded sign also allows inversion of a representative; this contributes only a factor of two per class.
\end{lemma}
\begin{proof}
For each length vector the description count is at most
\(\exp(h(\sum_iL_i-u)+aE(r\ell))\), up to polynomial factors. Multiply \eqref{fiber} for the $r$ classes and divide by $s(\ell)^r$. The additional term in the exponent is
\[
 -(h-D/2)\sum_i(\ell-L_i)\le0.
\]
A decoded representative determines its conjugacy class. Thus no enumeration of all geodesic spellings or all cyclic shifts is needed.
\end{proof}

\paragraph{Periodic descriptions.}
If an interval $p$ of a cyclic core is joined through an old diagram to an arc of $A^M$, put $P=M|A|_{\rm lit}$. After a cyclic shift of $A^{\pm1}$, write that arc as $T^jv$, with $0\le j\le M$ and $|v|_{\rm lit}<|T|_{\rm lit}=P/M$. For side labels $c,d$ and complementary core pieces $a,b$, the equality in $B$ is
\[
 R=a\,c\,T^jv\,d\,b.
\]
The six recorded group elements cost at most
\[
 |R|-|p|+|c|_B+|d|_B+2P/M.
\]
The saving is therefore
\begin{equation}\label{periodsave}
 |p|-|c|_B-|d|_B-2P/M.
\end{equation}
The outside pieces are recorded as group elements. The replacement may be longer than the deleted interval: its complete coding cost is already included in this saving.

\paragraph{Contacts between different cores.}
An equality $p=cq^{\pm1}d$ between intervals of two different cores has a seven-factor description with saving
\[
 \max(|p|,|q|)-(|c|_B+|d|_B).
\]
One records the shorter interval once, the four outside pieces, and the two connectors.

\paragraph{Self-contacts.}
Disjoint intervals are handled by the same repeated-factor idea. Partially overlapping intervals require more care. For
\[
 R=aXYZb,\qquad XY=c(YZ)d,
\]
the equation $(c^{-1}X)Y=Y(Zd)$ determines at most $M$ possibilities for $Y$ unless $c^{-1}X=1$. We use the finite cyclic centralizer theorem for nontrivial elements of $B$ \cite{Adian}. In the exceptional case $c=X$ and $d=Z^{-1}$ in $B$: this is the natural identification of the shared interval, not a rare event. For the reverse orientation,
\[
 XY=c(YZ)^{-1}d
 \quad\Longrightarrow\quad
 (Yd^{-1})^2=X^{-1}cZ^{-1}d^{-1},
\]
and $Y$ is uniquely determined because the square map is bijective in an odd-exponent group.

The cyclic case in which two arcs cover the whole core can be written
\[
 R=ABCD,\quad p=ABC,\quad q=CDA.
\]
For $p=e q f$, the two conjugacy equations are
\[
 A(BCf^{-1})A^{-1}=eCD,\qquad (e^{-1}AB)C=C(DAf).
\]
They allow deletion of $A$ or $C$, unless respectively
\[
 (e,f)=((CD)^{-1},BC)\quad\text{or}\quad(e,f)=(AB,(DA)^{-1}).
\]
The other two pieces can be solved for directly. For $p=e q^{-1}f$, deletion of either shared piece uses a square equation and has no such exception. Including the containment case gives the conservative bound
\begin{equation}\label{selfsave}
 \frac14\max(|p|,|q|)-(|c|_B+|d|_B)
\end{equation}
with at most six recorded group factors and multiplicity $O(M)$, unless there is a direct natural identification. These statements concern arcs of at most one traversal of a core.

For containment, If $p=XYZ$, $q=Y$, then the direct equation gives
\[
 (c^{-1}X)Y=Y(dZ^{-1}),
\]
with degeneracy exactly $c=X$, $d=Z$. The inverse equation gives
\[
 (YZd^{-1})^2=X^{-1}cZd^{-1}.
\]
For a full cyclic union $p=ABC$, $q=CDA$, the inverse equation $p=e q^{-1}f$ permits recovery of either shared piece through
\[
 (AK)^2=eK,\quad K=BCf^{-1}CD,\qquad
 (Cf^{-1})^2=B^{-1}A^{-1}eA^{-1}D^{-1}f^{-1}.
\]
Deleting the longest of the three or four pieces proves the conservative saving in \eqref{selfsave}. Unshared pieces can be recovered directly. The nondegenerate direct equations have at most $M$ solutions; the inverse equations have one. These are fixed multiplicities because $M$ is fixed before $\ell$ grows.

The direct natural identifications are kept as exceptions. In the diagram argument they give compatible cells, rather than exceptional random inputs; this is justified in Section~\ref{reductions}.

\section{Parameters and simultaneous generic exclusions}\label{parametersection}
The geometric perimeter parameter $n_0$ is independent of the actual exponent $M$. Condition A requires
\[
 |\partial\Pi|\ge n_0\rk(\Pi),
\]
so a power cell of perimeter $M\rk(\Pi)$ satisfies it whenever $M\ge n_0$. All geometric lemmas below use $n_0$; power relations and changes of exponents use $M$.

Fix $b\in(0,1)$, and choose $h_*$ with $D/2<h_*<H$. Choose an insertion hierarchy
\[
 n_0^{-1}\ll\zeta\ll\eps\ll\delta_g\ll\gamma\ll\beta\ll\alpha\ll1,
\]
including all the Condition-A inequalities. The symbol \(\ll\) means chosen sufficiently small in terms of the constants to its right. Here $\delta_g$ is a geometric parameter; write $h_g=\delta_g^{-1}$. Choose these constants from right to left, with $n_0$ last. In particular, require
\begin{equation}\label{hierarchy}
 3\gamma^{-1}\zeta<1,\quad
 2\eps^{-1}<\zeta n_0,\quad
 \zeta<\min\left\{\frac\eps{64},\frac{\eps\gamma}{32(1+\gamma)}\right\}.
\end{equation}
For the quantitative periodic comparisons, also require
\[
 \alpha<\frac1{100},\quad \beta<\frac\alpha{100},\quad
 \gamma<\frac\beta{100},\quad h_g\ge\frac{10^4}{\gamma},\quad
 \zeta<\frac{\eps\gamma^2}{1000}.
\]
Put
\[
 \begin{gathered}
 c_0=1-\beta,\qquad N_*=\eps n_0,\qquad
 T_0=\gamma^{-1}(2\zeta n_0+1),\\
 K_1=\gamma^{-1}+\alpha+2,\qquad K_0=2\gamma^{-1}+4+2\alpha.
 \end{gathered}
\]
Increase $n_0$ so that $\alpha<\zeta n_0$, $c_0(n_0-2)>2$, and
\begin{align}
 c_0N_*/2-2\alpha-\gamma^{-1}(2\alpha+1)&>h_g+2,\label{Q1}\\
 \frac{c_0\gamma}{2+\gamma}(c_0N_*/2-2\alpha)-K_1&\ge\gamma N_*/8,\label{Q2}\\
 K_0&<\gamma(\gamma N_*/8-5h_g/6-1)-1,\label{Q3}\\
 T_0&<N_*/2,\label{Q4}\\
 \frac{T_0(1+\gamma/2)}{c_0N_*/2-2\alpha}&<\gamma/2.\label{Q5}
\end{align}
These choices are compatible. The first three inequalities hold for sufficiently large $n_0$; the limiting inequalities in the last two follow from $\zeta<\eps\gamma^2/1000$. All are independent of the actual exponent $M$.

Fix a length-saving threshold
\begin{equation}\label{tauchoice}
 0<\tau<\min\left\{\frac{b}{16n_0},\frac{\eps b}{32},
                    \frac{\gamma\eps b}{8(1+\gamma)}\right\}.
\end{equation}
Next choose a tighter old hierarchy with
\[
 \beta_o\le\frac\gamma{8(1+\gamma)},
\]
and a coding hierarchy with
\begin{equation}\label{codingchoice}
 18\gamma_cD<\frac14h_*\tau,\qquad
 9\gamma_cD<\frac12(1-b)(h_*-D/2).
\end{equation}
Finally take the actual odd $M$ sufficiently large that all these hierarchies are valid for the old presentation and
\begin{equation}\label{Mchoice}
 \begin{gathered}
 h_M\ge h_*,\qquad H-h_M\le h_*\gamma/8,\qquad M>1000n_0/b,\\
 F_M:=(1-\beta_o)\eps M-2\zeta n_0
   >\max\{4+4\zeta n_0,\,8n_0/b\}.
 \end{gathered}
\end{equation}
The old periods are chosen using lengths and conjugacy, independently of these numerical hierarchies. Thus the same old presentation supports the tighter estimates. Thus all geometric and coding constants are fixed first, and then one exponent threshold $M_0$ is chosen. Every fixed odd $M\ge M_0$ is admissible. No parameter is subsequently sent to zero with $M$ fixed.

For each input choose a conjugacy-minimal word $R_i$, and put $L_i=|R_i|$. No special choice of its geodesic spelling is required. In the remainder of the proof set
\begin{equation}\label{insertionscale}
 d=\left\lfloor\frac{b\ell}{2n_0}\right\rfloor,
 \qquad \rho=d+\tfrac12.
\end{equation}
For sufficiently large $\ell$, if $L_i>b\ell$ then $\rho>2$ and $n_0\rho\le\min_iL_i$.

\begin{lemma}[Typical properties of the random relators]\label{goodevent}
With probability $1-\exp(-c\ell)$, for some $c>0$, all $L_i>b\ell$ and the following exclusions hold simultaneously for all cyclically shortest representatives of $g_i$ and $g_i^{-1}$, all cyclic positions, and intervals of at most one traversal.
\begin{enumerate}[label=(\roman*),leftmargin=*]
\item There is no periodic group description of the form used in \eqref{periodsave} with saving (reduction in total recorded length) at least $\tau\ell$.
\item No contact between different cores has saving at least $\tau\ell$ in the seven-factor description above.
\item No self-contact has saving at least $\tau\ell$ in \eqref{selfsave}, unless it is a direct natural identification.
\item Suppose an interval $p$ of $R_i$ is replaced, in $B$, by $cqd$, where $q$ is a specified segment of the infinite repetition of a freely reduced literal root $A$. There is no such description with
\begin{equation}\label{literalbad}
 h|p|-H|A|-h(|c|_B+|d|_B)\ge h_*\tau\ell.
\end{equation}
Here the root length, segment length and phase are $O(\ell)$.
\end{enumerate}
These are a finite list of fixed decoding schemes, not a quantifier over arbitrary encodings.
\end{lemma}
\begin{proof}
The length assertion is Lemma~\ref{cores}. For (i)--(iii), use Lemma~\ref{transfer} for the specified schemes. At most seven group elements are recorded, and at most two conjugacy classes are involved. The exponent is bounded above by
\[
 -h_*\tau\ell+18\gamma_cD\ell+7E(2\ell)+O(\log\ell),
\]
which is negative and linear in $\ell$ by \eqref{codingchoice}.

For (iv), record the four group elements consisting of the two complementary pieces and the two connectors, and record $A$ as a literal word. There are at most $C_m e^{H|A|}$ choices for $A$. Its phase and the length of $q$ determine $q$ and have polynomially many possibilities. The resulting code exponent is
\[
 h(L_i-|p|+|c|_B+|d|_B)+H|A|+4E(2\ell)+O(\log\ell).
\]
The conjugacy-fiber calculation therefore gives the same negative linear bound when \eqref{literalbad} holds. A decoded group element determines a conjugacy class, even when the encoding used a different shortest spelling. Thus no enumeration of geodesic spellings, ranks, strips, or diagrams is needed. A finite union handles the fixed tuple and the fixed list of schemes.
\end{proof}

\section{Adding the random relations at an intermediate rank and the diagram induction}\label{insertionsection}
Retain the old relations through rank $d$. At rank $\rho$ add the relators $R_1,\ldots,R_k$. Afterward, at each integer rank $j>d$, choose one representative of each conjugacy-up-to-inversion class of words of length $j$ which are simple in the preceding modified group, and impose its $M$th power. Denote these modified stages by $H_t$, where
\[
 t\in\{0,1,\ldots,d,\rho,d+1,d+2,\ldots\}.
\]
In particular $H_t=G_t$ for $t\le d$. A simple word is nontrivial and is not conjugate to a power of a shorter word or of an already selected period. Higher periods are reselected; we do not retain the original higher-rank relators as an additional family of cells.

The  remainder of Section \ref{insertionsection} and all of Section \ref{outersection} will be devoted to the proof of the following result. 
\begin{proposition}[Adding relations]\label{insertion}
Suppose the tuple of cyclic cores satisfies Lemma~\ref{goodevent}. With the preceding parameters and reduction convention below, every reduced disc or annular diagram over every $H_t$ is an A-map. A periodic boundary section is smooth of its root length if its root is simple in $H_t$, or is a selected period and the section has no compatible power cell. A cyclic periodic section must be a full power of its root. A marked new-relator section is smooth of rank $\rho$ if it has no compatible new cell. These assertions hold at every later rank with the same parameters.
\end{proposition}

\subsection{Reduction and marked sections}\label{reductions}
We first specify the surgeries used in the proof of Proposition~\ref{insertion}. Write $N_+(\Delta)$ for the number of positive cells. At a fixed stage $t$, the \emph{graded type} of $\Delta$ is the finite vector of its positive-cell counts, ordered by decreasing rank and compared lexicographically. Thus removing cells of rank $a$ and adding only cells of rank less than $a$ decreases graded type, but need not decrease $N_+(\Delta)$. This distinction will be respected throughout the two inductions.

\paragraph{Zero-refinement and labels.}
We use the zero-refinement of \cite[Section 11.5]{Olshanskii} and \cite[Section 2.2]{OS}. A thin band along an edge labelled $x$ consists of trivial cells with boundary labels $x1x^{-1}$; zero-labelled edges and cells separate the corners at vertices. Collapsing these bands recovers the original diagram. This separates coincident path occurrences and makes positive-cell closures embedded and disjoint, without changing any positive-cell label or rank. It also preserves the literal labels of the designated paths after identity letters are omitted. In particular, it neither creates nor duplicates positive cells. All boundary sections below retain their positions and orientations under this operation.

\paragraph{Compatible pairs and their replacement.}
For power cells of rank $j$, use the ordinary phase convention: oppositely oriented copies of $A^M$ are compatible when the relevant phase vertices are joined by a simple arc whose label is $1$ in the stage immediately below $j$. For new cells with signed cyclic boundary words $U,V$, based at the endpoints of a connecting arc labelled $T$, declare them compatible when
\begin{equation}\label{newcompat}
 UTVT^{-1}=1\quad\hbox{in }G_d.
\end{equation}
The arc is simple and its interior avoids the two cells. After zero-refinement, take a closed regular neighbourhood of the two cell discs and the arc. It is a disc: two disjoint discs joined by one embedded band have one boundary component and genus zero. Its oriented boundary has label $UTVT^{-1}$, with possible identity letters. A van Kampen disc over $G_d$ with this exact boundary label can therefore replace the neighbourhood. The replacement agrees with the rest of the diagram along the entire attaching boundary. It preserves the ambient disc or annular topology, every external boundary word, and every marked external position. It removes two rank-$\rho$ cells and adds only cells of rank at most $d$.

The same construction for a power pair uses the identity $A^MTA^{-M}T^{-1}=1$ below rank $j$, where $T=1$ there, and replaces the pair by cells of rank less than $j$. These are local disc replacements also in an annulus; the neighbourhood does not contain the annular hole. Neither replacement asserts a bound on the number of lower-rank cells introduced. Both strictly decrease graded type.

A new-relator boundary section carries a marking specifying its position in a full signed cyclic word $R_i$. Compatibility with a cell is tested using that full word in \eqref{newcompat}, with the orientation opposite to a cell attached outside the section. The witness must be an actual simple connecting arc whose interior avoids the cell and the section. The full-word marking is only part of this test; it does not attach another positive cell to the diagram. The following lemma explains why sections exposed from omitted cells are admissible.

A diagram is \emph{reduced} if no compatible pair exists after any zero-refinement. There is a reduced representative for each prescribed boundary datum and topology: choose one of minimal graded type. A compatible pair in a zero-refinement would admit the preceding replacement, contradicting minimality. Existence of a minimum follows from well-foundedness of the lexicographic order on a finite vector of nonnegative integers. We do not assert that every reduced diagram is minimal.

\begin{lemma}[Inheritance after excision]\label{inheritance}
Let $\Delta$ be reduced, and let $\Gamma$ be a subdiagram of a zero-refinement of $\Delta$.
\begin{enumerate}[label=(\roman*),leftmargin=*]
\item The diagram $\Gamma$ is reduced.
\item Suppose a boundary section of $\Gamma$ is copied from a cell $\Pi$ of $\Delta$ omitted from $\Gamma$. With its inherited full-relator marking or periodic phase, this section has no compatible cell in $\Gamma$.
\item A marked boundary section of $\Delta$ having no compatible cell retains that property on restriction to $\Gamma$.
\end{enumerate}
\end{lemma}
\begin{proof}
A zero-refinement of $\Gamma$ extends to one of $\Delta$ by adding zero-collars along its exposed boundary. Consequently a pair in any such refinement of $\Gamma$, including its simple connecting arc, is a pair in a refinement of $\Delta$. This proves (i). Notice that we do not project the arc to the unrefined diagram and claim that its projection is simple.

For (ii), let $\pi$ be a purported compatible cell in $\Gamma$. Its witness arc ends at the specified position on the exposed section. Continue it across the zero-collar to the corresponding position on $\partial\Pi$. The resulting arc has the same label and can be kept simple, with interior disjoint from $\pi$ and $\Pi$. The inherited marking makes the boundary compatibility equation exactly the pair equation for $\pi$ and $\Pi$; for a power cell it gives the same phase identity below its defining rank. The two cells are distinct because $\Pi$ is absent from $\Gamma$. This contradicts reducedness of $\Delta$. Part (iii) follows by the same extension of the witness arc, using the original boundary marking. This proves the inheritance assertions without any A-map or smoothness assumption; compare \cite[Lemma 3.17]{OS}.
\end{proof}

\paragraph{Natural self-contact identifications.}
The natural identifications from the self-contact count imply \eqref{newcompat}. For the strips used below, the connector equalities in these exceptional cases involve words of length at most $10\ell$: the side arcs have length at most $2\zeta\ell$ by \eqref{sides}, and each core piece has length at most $\ell$. Old stabilization puts these equalities in rank at most $\lceil30\ell/M\rceil$. By $M>1000n_0/b$ and \eqref{insertionscale}, this rank is at most $d$ for all sufficiently large $\ell$.

Here is the full-boundary calculation. If $R=PQ$ and $R'=QP$ are the two based readings, the natural connector is represented by $P$ or by $Q^{-1}$. In either case
\[
 RP=P R',\qquad RQ^{-1}=Q^{-1}R'
\]
in the free group, and hence in $G_d$. The exceptional connector equalities just established in $G_d$ give the same identities with the actual connector in place of $P$ or $Q^{-1}$. Read the other cell with the opposite sign to obtain \eqref{newcompat}. This covers both prefix transport and transport through the complementary suffix. In particular, if
\[
 U=ABCD,\quad W=CDAB,\quad T=(CD)^{-1},
\]
then $UT=TW$ in the free group. Thus this full-cover exception gives a pair with oppositely oriented $W$, although $T$ need not equal $AB$ in $G_d$. The connector is an actual side arc of the exposed strip, so the required simple witness is present; we need not turn an arbitrary group equality into a simple path in the ambient diagram.

For a self-strip on one cell, the two physical occurrences are disjoint intervals of the same oriented cell contour. Reading both forward in that single signed core changes the strip relation into the inverse self-contact equation. Its counting argument has no direct natural-identification exception. Thus a self-strip is handled by that equation, not by treating a cell as a distinct-cell pair with itself.
\subsection{Shortcut discs and complementary arcs}\label{shortcutcuts}
The following form of shortcut excision records the topological input used in the cell-count induction. We give a plane-graph proof, including the annular case. In particular, we do not assume that an intermediate return of a shortcut cuts off a contractible loop.

\begin{lemma}[Shortcut excision]\label{shortcutexcision}
Let $\Delta$ be a reduced disc or annular diagram with $N$ positive cells. Suppose an interval $u\subset\partial\Pi$ of a positive cell has a shorter path $v$ homotopic to $u$ relative to its endpoints. After zero-refining and replacing $v$ by a path of no greater length if necessary, while retaining $u$ and its label, one of the following disc regions is available (with $v$ denoting the replacement):
\begin{enumerate}[label=(\roman*),leftmargin=*]
\item a disc $E$ with contour $uv^{-1}$ which omits $\Pi$;
\item a disc $E'$ with contour $vu'$, where $\partial\Pi=uu'$ and $E'$ omits $\Pi$.
\end{enumerate}
The available disc contains at most $N-1$ positive cells and is reduced. Its cell-contour side carries the inherited marking and has no compatible cell. If instead $u$ lies on a designated boundary section of $\Delta$, a disc with contour $uv^{-1}$ is obtained in the same way and has at most $N$ positive cells. In both assertions every positive cell of the resulting disc is an original cell, used at most once.
\end{lemma}
\begin{proof}
First make the closure of $\Pi$ embedded. If $u$ is an entire contour, insert a zero edge at its basepoint and leave that edge outside $u$, so that its endpoint occurrences are distinct. Now give $\Pi$ a collar of trivial cells, as in Section~\ref{reductions}. Take $\partial\Pi$ on the inner bank of the collar and keep the given competitor on the old one-skeleton, joining its terminal occurrences to the endpoints of $u$ by zero edges. Thus $v$ meets $\partial\Pi$ only at those endpoints; it may still have self-intersections. Labels, lengths and the relative homotopy class are unchanged, including under the basepoint convention just described.

Let $K$ be the finite connected plane graph consisting of the edges traversed by $u$ and $v$, with their multiplicities forgotten. Every interior vertex of $u$ has degree two in $K$, and the closed walk $C=uv^{-1}$ traverses each edge of $u$ exactly once. Cap the boundary components of the ambient surface by discs. We obtain a sphere with one distinguished cap in the disc case and two in the annular case. The winding index of $C$ is constant on each face of $K$; normalize it to be zero on the first cap. In the annular case it is also zero on the other cap, since $C$ is null-homotopic in the annulus. The indices on the two faces incident to $u$ differ by one. These are the same two faces all along $u$, because its interior vertices have degree two in $K$. Consequently at least one such face $F$ has nonzero index and meets no cap. This is the step that excludes the annular hole.

Orient the facial boundary walk of $F$ so that it has the form $uw^{-1}$. We claim that
\[
 |w|\le |v|<|u|.
\]
Indeed, in a connected plane graph a facial boundary walk uses any non-bridge edge at most once, and any bridge at most twice. Every edge of $K$ occurs in $C$, and a closed walk using a bridge must traverse it at least twice. Thus, edge by edge, the number of occurrences in the facial boundary is at most the number in $C$. Each edge of $u$ occurs once in both and never in the given $v$. Removing these occurrences proves the displayed inequality, also when some edges have length zero. Notice that $w$ need not be a subpath of the given $v$.

The face $F$ is a face of a connected graph in a sphere and contains no cap. Separate the two banks of each bridge and split repeated boundary vertices by zero-bands and zero-corners. Its closure is thereby realized as an embedded disc $D$ with contour $uw^{-1}$. These operations only separate occurrences of edges and vertices; each positive cell in $D$ is a cell of $\Delta$ occurring once. In particular, this is an exposed subdiagram of a zero-refinement, not a filling obtained by lifting to a covering space. Replace $v$ by $w$.

Since $K\cap\partial\Pi=u$, the interior of $\Pi$ lies entirely on one side of this facial boundary. If $D$ omits $\Pi$, take $E=D$, giving (i). Otherwise $\Pi$ is contained in $D$ and meets its contour exactly along $u$. Removing $\Pi$ through that boundary interval leaves a disc $E'$ with contour $vu'$, where $\partial\Pi=uu'$, giving (ii). Zero-collars retain a disc also when the remainder has no positive cells. In either case the resulting disc omits $\Pi$, so it has at most $N-1$ positive cells. Reducedness and the assertion about the inherited marking follow from Lemma~\ref{inheritance}.

For a designated boundary interval, attach a trivial collar to the ambient boundary, put $u$ on its new boundary bank, and keep the given competitor on the old graph. The same argument produces $D$ with contour $uv^{-1}$ and at most $N$ positive cells. Here no distinguished positive cell is to be removed. The zero-edge convention above handles an entire cyclic contour, and Lemma~\ref{inheritance} again gives reducedness and the inheritance of any no-compatible-cell hypothesis.
\end{proof}

\subsection{Removal at a periodic boundary}\label{boundaryremovalsection}
We will also need a surgery which decreases graded type. Unlike shortcut excision, it is used only after the cell-count induction at the current stage has been completed.

\begin{lemma}[Boundary power-cell removal]\label{boundaryremoval}
Let a disc or annular diagram over $H_t$ have a designated section $q$ labelled $B^s$, where $B$ is a selected period of rank $j\le t$. Then it can be replaced by a reduced diagram with the same topology and the same other designated boundary words, in which this section is labelled $B^{s'}$, where $s'\equiv s\pmod M$, and has no compatible $B$-cell. If a compatible cell is removed, the graded type strictly decreases. All other boundary labels and marked endpoints are preserved, but no bound on the total positive-cell count is asserted.
\end{lemma}
\begin{proof}
Suppose $\Pi$ is compatible with $q$. Let $t_0$ be its simple phase arc, with label $T=1$ in $H_{j^-}$, the stage immediately preceding $j$. Write $s=s_1+s_2$ at its endpoint on $q$. Slit along $t_0$, go around $\Pi$, and return along the other bank of the slit. Removing $\Pi$ changes the designated boundary path to one with label
\[
 W=B^{s_1}T B^{\eta M}T^{-1}B^{s_2},\qquad \eta\in\{1,-1\},
\]
where the sign is determined by the orientation of the exposed cell. The rest of the boundary is unchanged. In $H_{j^-}$ we have
\begin{equation}\label{boundaryreplacement}
 W=B^{s+\eta M},
\end{equation}
because $T=1$ there. Notice that this equality does not use $B^M=1$ below rank $j$.

Glue a diagram of \eqref{boundaryreplacement} along the exposed path $W$. It contains only cells of rank less than $j$ and replaces that path by $B^{s+\eta M}$. Topologically the cell and its phase arc were removed through a notch at the designated boundary, and the equality diagram closes that notch. Thus the original disc or annular topology is retained. If $q$ is an entire contour, perform the construction with a basepoint slit and then identify its two banks. If the new exponent is zero, retain a zero-edge boundary contour rather than cap an annular hole. This convention also retains the designated boundary partition in the zero-length case.

The operation removes one rank-$j$ cell and adds only lower-rank cells, so it strictly decreases graded type. Reduce compatible pairs using the replacements of Section~\ref{reductions}, and repeat whenever the new designated section has a compatible $B$-cell. Every nontrivial step decreases the finite graded type, so this process terminates. It leaves the other boundary words fixed and changes the exponent only by multiples of $M$. This is the boundary surgery of \cite[Lemma 3.18]{OS}, written here with trivial zero-subgroups. It makes no assumption about the labels of other positive cells.
\end{proof}

We will apply this lemma separately to each auxiliary diagram in Section~\ref{outersection}. Admissibility in one annulus does not establish admissibility in a subsequently constructed disc. Nor does the lemma assert that another marked section retains its no-compatible-cell property after lower-rank fillings have been added; that property must be checked when it is required.

\subsection{Strict ordering of the diagram induction}
Fix a rank $t$, and assume the periodic comparison and period-recognition assertions at every strictly smaller rank. At $t=\rho$ these are old results. At this fixed rank induct on the number $N$ of positive cells, proving, in order:
\begin{description}[leftmargin=1.5cm,style=nextline]
\item[$D_N$] Every reduced disc or annular diagram with at most $N$ cells is an A-map.
\item[$S_N$] In such a diagram the following sections are smooth: a marked $R_i$-interval with no compatible new cell; a periodic section rooted in a selected period with no compatible power cell; and a periodic section rooted in a word simple in $H_t$.
\end{description}
For cyclic sections in the second and third classes the label must be a full power $A^s$. In particular, \emph{cyclic powers of earlier periods with no compatible cell are included}. They are required after compatible-cell removal in the outer annular arguments. A full cyclic $R_i$ is allowed in the first class, with its compatibility condition. Their respective smoothness ranks are $\rho$ and $|A|$.

All uses of smoothness in proving $D_N$ occur in proper subdiagrams, so use $S_{N-1}$. After $D_N$ has been proved, local geodesicity in $S_N$ may use $D_N$, but its contiguity tests still use proper strips and $S_{N-1}$. For $N=0$ the claims follow from free diagram geometry. Every selected period and every current-simple word is cyclically freely reduced, and the same is true of each cyclic core.

Inside this induction, we only expose zero-refined subdiagrams and use Lemma~\ref{inheritance}; we do not first replace them by arbitrary reduced fillings. The pair and boundary replacements may increase $N_+$, so they are used to establish the reduction convention or, after $D_N,S_N$ hold for all $N$, in the outer algebraic arguments. In particular Lemma~\ref{boundaryremoval} is not used to claim that a new diagram has fewer than $N$ cells.

\subsection{Every relevant mixed strip is old}
Let $\Gamma$ be a contiguity strip tested in $D_N$ or $S_N$, with contour $cpdq$. By definition it is a disc subdiagram omitting its source cell and, in a cell-to-cell test, its target cell. Thus $N_+(\Gamma)\le N-1$, and $N_+(\Gamma)\le N-2$ when the two endpoint cells are distinct. A self-contiguity omits its single endpoint cell and still satisfies the first bound. These statements concern the strip itself; the union of a self-strip with its cell can surround an annular hole.

We justify using the strip as its own ambient map. Every defining bond lies in $\Gamma$. For a zero-bond its zero-band, positive end edges and zero side paths are unchanged when an omitted endpoint cell is replaced by an exposed boundary section. For a positive bond, its principal cell and its two recursively defining lower bonds or strips are all contained in that bond. Induction on this recursive construction preserves the disjointness, the degree inequalities and every comparison between internal principal-cell ranks. A rank inequality involving an endpoint cell is retained if that endpoint is still a cell; if the endpoint becomes a boundary section, that inequality is no longer required by the bond definition. Boundary smoothness ranks are not additional inequalities in the definition of a bond. Zero-collars resolve any boundary coincidences without changing these data. Hence the original bonds still define $\Gamma$ as a boundary-to-boundary contiguity strip in its own map; see \cite[Section 2.3]{OS}.

By Lemma~\ref{inheritance}, $\Gamma$ is reduced and every section inherited from an omitted cell has no compatible cell in $\Gamma$. A no-compatible-cell hypothesis on an ambient boundary section is likewise inherited. If the target is rooted in a word simple in $H_t$, we keep that same ambient stage $t$ while applying $S_{N-1}$: it is a statement about diagrams over $H_t$ with fewer cells, so no new simplicity assertion at $\rk\Gamma$ is needed here. It follows from $D_{N-1}$ and $S_{N-1}$ that $\Gamma$ is an A-map and its two indicated sections are smooth of ranks $a,b$.

Now apply the boundary-to-boundary case of \cite[Lemma 2.1]{OS} to this A-map and its just-verified contiguity structure, in both directions. That case permits the strip to be the whole ambient disc; properness was needed relative to the original diagram to invoke the induction. For its side labels $c,d$ this gives
\begin{equation}\label{sides}
 |c|,|d|<\zeta n_0\min(a,b).
\end{equation}
The rank-drop lemma \cite[Lemma 2.14]{OS} now gives
\begin{equation}\label{striprank}
 \rk(\Gamma)<3\gamma^{-1}\zeta\min(a,b)<\min(a,b).
\end{equation}
If a section comes from a new relator, $\rk(\Gamma)<\rho$, so $\Gamma$ is an old diagram. No conclusion about the full ambient mixed diagram was used to obtain this fact.

In particular, \eqref{sides} retains the rank of the source. This dependence will pay for both connectors in the power-to-new-relator case.

\subsection{The three mixed cases and their budgets}
\paragraph{Power source, new target.}
Let the source have label $A^M$, rank $j=|A|$, and degree at least $\eps$. Let $p$ be its contiguity arc, $q$ the new-relator arc, and $C=|c|+|d|$. By \eqref{striprank}, the strip is old and lies below both $j$ and $\rho$. The word $A$ is simple in the relevant old stage: below $j$ for a retained old period; in $G_d$ for a newly selected period, since a failure of simplicity there would persist in the modified stage where it was selected.

Old periodic quasigeodesicity, with the tighter old parameter, gives
\[
 (1-\beta_o)\eps Mj\le (1-\beta_o)|p|\le |q|+C,
 \qquad C<2\zeta n_0j.
\]
Consequently
\begin{equation}\label{repairedbudget}
 j\le\frac{|q|}{F_M},\qquad
 2j+C\le\frac{2+2\zeta n_0}{F_M}|q|<\tfrac12|q|.
\end{equation}
The root, its terminal remainder, and both connectors therefore cost less than $|q|/2$ as group elements. If $|q|\ge(1+\gamma)\rho$, the periodic description saves more than $\rho/2>\tau\ell$ for large $\ell$, contrary to Lemma~\ref{goodevent}(i). Also $j\le\ell/F_M<\rho$, so the source is itself old. Both conclusions follow without an ambient A-map assumption.

\paragraph{New source, periodic target.}
Let $p$ be the new-relator arc, $|p|\ge\eps L$, and let $q$ be periodic with root $A$ of length $a$. The strip is old. The same old-stage simplicity argument applies to $A$, and hence
\begin{equation}\label{excess}
 |p|\ge(1-\beta_o)|q|-C,
 \qquad C<2\zeta n_0\min(\rho,a)\le2\zeta L.
\end{equation}
Since the core is cyclically $B$-geodesic, also $|q|\ge|p|-C\ge(\eps-2\zeta)L$. If the smoothness bound fails, then $|q|\ge(1+\gamma)a$. Recording the root \emph{literally} has entropy cost $Ha$, so the entropy saving is at least
\begin{align}
 h|p|-Ha-hC
 &\ge\left[h(1-\beta_o)-\frac H{1+\gamma}\right]|q|-2hC\label{paid}\\
 &\ge\frac{h\gamma}{2(1+\gamma)}|q|
 \ge\frac{h\gamma\eps}{4(1+\gamma)}L
 >h_*\tau\ell.\nonumber
\end{align}
Here \eqref{hierarchy}--\eqref{Mchoice} give
\(h(1-\beta_o)-H/(1+\gamma)\ge3h\gamma/[4(1+\gamma)]\) and
\(C\le\gamma|q|/[8(1+\gamma)]\).
Moreover $|q|\le(L+C)/(1-\beta_o)=O(\ell)$, and $a\le|q|$, so the polynomial auxiliary bound in Lemma~\ref{goodevent}(iv) applies. This contradicts that lemma. Equation~\eqref{paid} explicitly pays for both the replacement excess and the growth gap; using the coarser insertion parameter $\beta$ here would fail.

\paragraph{Contiguity between sections of the new relators.}
One arc has length at least $\eps L\ge\eps b\ell$, whereas $C<2\zeta n_0\rho\le2\zeta L$. An independent contact, or an essential self-contact, therefore has saving at least $\eps b\ell/8>\tau\ell$. Lemma~\ref{goodevent}(ii)--(iii) excludes it. The remaining direct natural identification gives an old compatible pair by \eqref{newcompat}, contrary to reduction or boundary admissibility. A self-contiguity of one new cell gives the inverse case, which has no exceptional identification.

\subsection{Proof of \texorpdfstring{$D_N$ and $S_N$}{DN and SN}}
Condition A1 follows from $M\ge n_0$ and the insertion scale. For A2, put $r=\rk\Pi$, take a cell-boundary arc $u$ of length $s\le\max(r,2)$, and suppose it has a shorter homotopic competitor $v$. Apply Lemma~\ref{shortcutexcision}.

In its first case the shortcut disc $E$ omits $\Pi$, so $D_{N-1}$ makes it an A-map. If $s\le r$, then $|\partial E|=s+|v|<2r$. For any positive cell $\pi$ in $E$, the perimeter inequality and A1 give
\[
 c_0n_0\rk\pi\le c_0|\partial\pi|<|\partial E|<2r.
\]
Thus $\rk E<2r/(c_0n_0)<r$. For a power cell, the label of $u$ is an interval of one cyclic period. Replacing it by the shorter label of $v$ gives, below the defining rank, a shorter word conjugate to that period, contrary to its simplicity. For a new cell $r=\rho$, and $E$ is old; the equality $u=v$ therefore holds in $B$, contrary to cyclic $B$-geodesicity. The only extra test allowed by $\max(r,2)$ has $s=2$ and $r=1$. Its shortcut perimeter is at most three, so the same perimeter inequality and $c_0n_0>3$ rule out every positive cell. Equality in the free group then contradicts cyclic free reduction.

In the second case of Lemma~\ref{shortcutexcision}, the disc $E'$ has contour $vu'$ with $\partial\Pi=uu'$ and at most $N-1$ positive cells. The \emph{complementary} arc $u'$ carries the full-relator marking or the phase inherited from $\Pi$. It is admissible by Lemma~\ref{inheritance}(ii); any compatible cell would pair with the omitted $\Pi$. Consequently $D_{N-1}$ and $S_{N-1}$ apply to $E'$ and $u'$. The smooth-section inequality \cite[Lemma 2.7]{OS} gives
\[
 c_0|u'|\le|v|<s.
\]
On the other hand, A1, $s\le2r$ and the chosen inequality $c_0(n_0-2)>2$ imply
\[
 c_0|u'|=c_0(|\partial\Pi|-s)
 \ge c_0(n_0r-s)\ge c_0(n_0-2)r>2r\ge s,
\]
a contradiction. This uses no A2 or A3 assertion for the current $N$-cell diagram. The annular case is covered by the contractible-region assertion in Lemma~\ref{shortcutexcision}.

For A3, the tested strip is proper. If one cell is new, the three preceding calculations give the required bound. If both cells are powers, \eqref{striprank} lowers the strip below both period ranks. The earlier-rank periodic comparison lemma gives either the required target-arc bound or period identification and compatibility. Distinct compatible cells contradict reduction. For a self-contiguity, both arcs inherit the same periodic orientation from the one cell; the orientation clause of that comparison excludes this case. One cannot use the distinct-cell pair argument for a self-contiguity. This proves $D_N$.

For local geodesicity in $S_N$, use the boundary case of Lemma~\ref{shortcutexcision}. Its disc has at most $N$ cells and is reduced, so the already established $D_N$ applies. A shortcut on a marked $R$-arc of length at most $\rho$ gives perimeter less than $2\rho$; the preceding perimeter calculation makes this disc old, contradicting cyclic $B$-geodesicity. For a selected period of length $j$, the same calculation lowers a shortcut on at most one period below $j$ and contradicts simplicity before selection. For a word simple in $H_t$, replacing an interval of one cyclic period directly gives a shorter conjugate in $H_t$. The additional length-two test when the root has length one has perimeter at most three and reduces to free geometry. This is the precise point at which $D_N$, rather than merely $D_{N-1}$, is required.

Finally, a contiguity strip tested in $S_N$ omits its source cell. Apply $D_{N-1}$ and $S_{N-1}$ to that strip, then \eqref{striprank}. The three mixed cases apply without change. For power--power contacts use the earlier-rank comparison: an existing-period target would have a compatible cell, and a current-simple target would be conjugate to a selected period, both excluded by the hypotheses. This includes cyclic full powers of earlier periods when no compatible cell is present; the periodic label continues across the chosen basepoint. The same argument handles cyclic simple powers and marked cyclic $R$-contours. Thus $S_N$ follows, completing the inner induction.

\section{Periodic comparison and the outer induction}\label{outersection}
Fix a stage $t$. Assume the comparison assertions at all smaller stages and the A-map and smoothness assertions established in Section~\ref{insertionsection} at $t$, for all cell counts. We prove the comparison assertions at $t$. Only the geometric lemmas from Section~\ref{inputs} are invoked as statements valid independently of the presentation. The group comparisons below extend the arguments of \cite[Sections 18--19]{Olshanskii} to the inserted presentation.

Put $c_0=1-\beta$ and $\sigma=1/2+\alpha$. Compatibility of two oppositely read $A$-periodic sections means that two corresponding phase vertices are joined by a path representing $1$ in $H_t$. When $A$ is a selected period of rank $j$, cell reduction requires this identity below $j$. Our uses of that reduction follow rank decrease to a stage below $j$.

\subsection{Root descent, torsion, and boundary surgery}
Descent in word length writes every nontrivial element as a conjugate of $A^s$, where $A$ is simple in $H_t$ or a selected period. If a current-simple word $A$ had finite order, a reduced disc for $A^s=1$, $s\ne0$, would have a smooth cyclic contour. The smooth-length inequality would give $c_0|A^s|\le0$, a contradiction. Every finite-order element is therefore conjugate to a power of a selected period. Its order divides $M$, so $H_t$ has no involutions.

If current-simple words $A,B$ satisfy $A\sim B^s$, simplicity gives $|A|\le|B|$, whereas smooth annular comparison gives $|A|\ge c_0|s||B|$. Since $c_0>1/2$, necessarily $s=\pm1$. Simplicity in both directions now gives $|A|=|B|$.

For nontrivial conjugate words $X,Y$, the annular cut (G6), followed by contour arcs of length at most $|X|/2$ and $|Y|/2$, gives a conjugator of length at most $(1/2+\gamma)(|X|+|Y|)<\sigma(|X|+|Y|)$. An essential loop of length zero would make the conjugacy class trivial, so (G6) applies.

For a selected $B$-periodic boundary apply Lemma~\ref{boundaryremoval}, with the other boundary words and their marked endpoints fixed. Its formula \eqref{boundaryreplacement} shows exactly why the exponent changes by a multiple of $M$ and why the filling lies below $|B|$. At this point $D_N,S_N$ have been established for every $N$ at stage $t$, so possible growth in the number of lower-rank cells is harmless. The reduced terminal diagram is an A-map and its nonempty $B$-section is smooth, since it has no compatible cell. We perform this removal afresh in every auxiliary diagram where it is needed. If the power side is empty, retain its zero-length boundary and apply length comparison to the other smooth side.

\subsection{Three simultaneous comparisons}
\begin{lemma}\label{pcm}
For reduced disc diagrams over $H_t$ the following statements hold. All roots mentioned are simple in $H_t$.
\begin{enumerate}[label=(\alph*),leftmargin=*]
\item\label{pcm:p} If the contour is $p_1q_1p_2q_2$, the labels of $q_1$ and $q_2^{-1}$ are $A$-periodic, and
\[
 |p_1|,|p_2|<\alpha|A|,\qquad
 |q_1|,|q_2|>(5h_g/6+1)|A|,
\]
then the two periodic sections have compatible phase vertices.
\item\label{pcm:c} Suppose $u,v$ are positive integers, $r=\min(u,v)$, and
\begin{equation}\label{centralizerbudget}
 Z_1A^uZ_2=A^v,\qquad
 |Z_1|+|Z_2|<[\gamma(r-5h_g/6-1)-1]|A|.
\end{equation}
Then $Z_1$ and $Z_2$ represent powers of $A$.
\item\label{pcm:m} Suppose $a=|A|\ge b=|B|$, the contour is $p_1q_1p_2q_2$, and $q_1,q_2$ are respectively $A$- and $B$-periodic, with
\[
 |p_1|,|p_2|<\alpha b,\qquad |q_1|>3h_ga/4,\quad |q_2|>h_gb.
\]
Then $A\sim B^{\pm1}$ and $a=b$. If $q_1$ starts with $A$ and $q_2^{-1}$ starts with $B^{-1}$, and $P$ labels $p_1$, then $A=P^{-1}B^{\pm1}P$ in $H_t$.
\end{enumerate}
\end{lemma}
\begin{proof}
At fixed $t$ induct on the sum $J$ of root lengths, in the order (a), (b), (c) at each $J$. For (a) and (b), $J=2|A|$. Thus (c) may use (b) at $2|B|=J$ when the roots have equal lengths.

\emph{Same-period compatibility.}
Put $a=|A|$. Align the initial phases by shifting an endpoint by at most $a/2$, replacing $A$ by a cyclic shift when needed. The resulting connector labels $D,E$ have $|D|<\sigma a$, $|E|<\alpha a$, and both periodic sides remain longer than $5h_ga/6$. If $D=1$ in $H_t$, compatibility holds.

Suppose $D\ne1$. Interchanging the sides if necessary, their labels can be written $QV$ and $Q^{-1}$, with $s=|Q|$ the shorter length. Glue the matching $Q$-paths to obtain an annulus with contours labelled $D$ and $VE$, up to orientation, and reduce it. Its core represents a nontrivial element. Smoothness before gluing gives
\[
 |V|\le(c_0^{-1}-1)s+c_0^{-1}(|D|+|E|).
\]
Apply (G6), and extend the cut by at most half of the $D$-contour and half of the $E$-path to connect the initial vertex of $Q$ to a vertex on $V$. Its label $S'$ satisfies
\begin{equation}\label{annularshort}
 |S'|<\gamma(c_0^{-1}-1)s+2a.
\end{equation}
The path $Q$ followed by the appropriate initial part of $V$ is an $A$-periodic word $W$, with $|W|\ge s$. The two paths to the same endpoint differ by a power of the annular core, so
\begin{equation}\label{annularequation}
 W=D^\nu S'\quad\text{in }H_t
\end{equation}
for some integer $\nu$. Reductions preserve this group equation, which is the datum used below.

For clarity, the old $Q$-seam need not remain a literal path after pair reduction. Each replacement is a disc with the same attaching word. A path crossing that disc can be continued through the replacement with the same endpoints and the same label as a group element of $H_t$: the old and new subpaths are each equal in $H_t$ to a path on their common attaching boundary. The outer boundary words, including the chosen endpoint on $V$, stay fixed. Carry the old seam followed by the relevant prefix of $V$ through these replacements. Its label still represents $W$, although no length estimate is claimed for this carried path. In the resulting annulus it and the short path labelled $S'$ differ by an integer winding about the based $D$-contour, because the fundamental group of the underlying annulus is cyclic. Reading that loop gives exactly \eqref{annularequation}. Only $S'$ is assigned the new geometric length bound. This argument uses the disc replacements of Section~\ref{reductions} after completion of the inner induction, so no cell-count bound is needed here.

Write $D=YB^rY^{-1}$ by root descent. If $B$ is a selected period, remove compatible $B$-cells from this conjugacy annulus. Its power contour is then smooth. Since $D\ne1$, that contour is nonempty, and
\[
 b=|B|\le|B^r|\le c_0^{-1}|D|<c_0^{-1}\sigma a<3a/5.
\]
The short-conjugator estimate permits $|Y|<a$. Equation~\eqref{annularequation} gives a disc with $A$-periodic side $W$, a full $B$-power side, and total connector length less than
\[
 \gamma(c_0^{-1}-1)s+4a.
\]
If $B$ is a selected period, remove compatible $B$-cells again in this new disc, and reduce. This preserves $W$ and the connectors, changing only the $B$-exponent. The $B$-side has length at least $c_0|W|$ minus the connector total, so it remains nonempty and long enough for (G5). Reverse the sign of the root if necessary.

Trim with minimum smooth rank $b$. Each long side loses less than
\[
 (c_0^{-1}-1)s+5\gamma^{-1}a.
\]
The retained $A$-side is longer than
\[
 (2-c_0^{-1})s-5\gamma^{-1}a>3h_ga/4.
\]
The retained $B$-side is longer than $h_gb$, by (G2), $b<3a/5$, and the parameter inequalities. Crucially, both new connectors are shorter than $\alpha b$. If $B$ is current-simple, apply (c) at $a+b<2a$. If $B$ is a selected period, (G4) first lowers the trimmed diagram below its rank $b$, where both roots are simple and the earlier-rank comparison applies. Either alternative makes $A$ conjugate to $B^{\pm1}$, contradicting simplicity and $b<a$. Thus $D=1$ and (a) follows.

\emph{Long power equations.}
Take a reduced disc for \eqref{centralizerbudget}. Its two periodic sides are smooth. Trimming loses less than
\[
 \gamma^{-1}(|Z_1|+|Z_2|+a)<(r-5h_g/6-1)a
\]
from each. The retained sides meet the thresholds in (a). Reading the original connector through the compatible phase vertices yields $Z_1=A^e$; the original equation then gives $Z_2=A^{v-u-e}$. This proves (b).

\emph{Different periodic roots.}
Choose cyclic shifts of the roots without moving the endpoints, and write the boundary equation as
\[
 P A^l A_0 Q B_0^{-1}B^r=1,
\]
where $|A_0|<a$, $|B_0|<b$, and $r>h_g-1$. Cyclic shifting and substitution give
\begin{equation}\label{mixedidentity}
 (PA^{-1}P^{-1})B^{-r}
 (B_0Q^{-1}A_0^{-1}AA_0QB_0^{-1})B^r=1.
\end{equation}
The two connector factors have total length less than $4a+2b+4\alpha b$. The original length comparison gives $3h_ga/4<c_0^{-1}(r+2)b$, hence
\[
 \frac{4a+2b+4\alpha b}{b}
 <\frac{16(r+2)}{3c_0h_g}+2+4\alpha
 <\gamma(r-5h_g/6-1)-1.
\]
The last inequality holds for every $r>h_g-1$ under $h_g\ge10^4/\gamma$: its left side is less than $14r/h_g$, while the final expression is greater than $\gamma r/8$. Apply (b) with root $B^{-1}$ to \eqref{mixedidentity}. It gives $PA^{-1}P^{-1}=B^e$. Simple-power conjugacy forces $e=\pm1$ and $a=b$, and gives the stated based conclusion. This completes the simultaneous induction.
\end{proof}

\subsection{Orientation}
\begin{lemma}\label{orientation}
If both $q_1$ and $q_2$ in the contour $p_1q_1p_2q_2$ are periodic in the same current-simple root $A$, and $|p_1|,|p_2|<\alpha|A|$, then
\[
 \max(|q_1|,|q_2|)\le h_g|A|.
\]
\end{lemma}
\begin{proof}
Suppose the larger length exceeds $h_ga$, where $a=|A|$. Length comparison makes the other exceed $c_0h_ga-2\alpha a>3h_ga/4$. Apply Lemma~\ref{pcm}(c), interchanging the sides to meet its thresholds. Choose $q_1$ to start with $A$ and $q_2^{-1}$ to start with $C=V^{-1}A^{-1}V$, where $|V|<a$. Put $P=\operatorname{Lab}(p_1)$ and $Z=VP$. The based conclusion gives
\[
 A=Z^{-1}A^{\pm1}Z,\qquad |Z|<(1+\alpha)a.
\]
If $Z$ centralizes $A$, write $q_1=A^rA_0$ and $q_2^{-1}=C^sC_0$. The boundary equation gives
\[
 A^{r+s}=VC_0\operatorname{Lab}(p_2)^{-1}A_0^{-1}Z^{-1}.
\]
Its right side has length less than $(4+2\alpha)a$, whereas $(r+s)a>(h_g-2)a$. This contradicts quasigeodesicity of simple powers and $c_0(h_g-2)>4+2\alpha$.

Otherwise $Z^{-1}AZ=A^{-1}$. Apply Lemma~\ref{pcm}(b) to $Z^2A^rZ^{-2}=A^r$ for a sufficiently large positive integer $r$. It gives $Z^2=A^e$. Conjugation by $Z$ gives $A^e=A^{-e}$. Since $A$ has infinite order, $e=0$. Now $Z^2=1$, and absence of involutions gives $Z=1$, contradicting $A=A^{-1}$.
\end{proof}

\subsection{One period compared with a very long periodic side}
Recall $N_*=\eps n_0$ and the constants $K_0,K_1,T_0$ from Section~\ref{parametersection}.
\begin{lemma}\label{longshort}
Let $A,B$ be current-simple, of lengths $a,b$. Suppose a reduced disc has contour $p_1q_1p_2q_2$, where $q_1,q_2$ are respectively $A$- and $B$-periodic, and
\[
 |p_1|,|p_2|<\alpha b,\qquad
 |q_1|>(1+\gamma/2)a,\qquad |q_2|>N_*b/2.
\]
Then $A\sim B^{\pm1}$. If $q_1$ starts with $A$ and $q_2^{-1}$ starts with $B^{-1}$, the conjugator is the label of $p_1$, as in Lemma~\ref{pcm}(c).
\end{lemma}
\begin{proof}
First suppose $a<b$. Length comparison gives $|q_1|>(c_0N_*/2-2\alpha)b$. Trim with minimum rank $a$, losing at most $\gamma^{-1}(2\alpha+1)b$ from each side. By \eqref{Q1} the retained sides are long enough for Lemma~\ref{pcm}(c) with the roots interchanged, and the new connectors have length less than $\alpha a$. That lemma forces $a=b$, a contradiction. Hence $a\ge b$.

Cyclically shift $A$ so that $q_1$ begins with $A$, and $B$ so that $q_2$ ends with $B$. These changes are unnecessary under the hypotheses for the based conclusion. Write $q_1=uv=wu'$, where $|u|=|u'|=a$. Its periodicity gives literal equality $V=W$ for the labels of $v,w$. Both have length $L_0=|q_1|-a$. By the interior-point cutting lemma, the endpoints of $u$ and $w$ can each be joined to $q_2$ by a path shorter than $\gamma^{-1}b$. These are separate cuts and need not be disjoint from each other.

Choose the $B$-phase at the terminal endpoint of $q_2$. Absorb its initial partial period into $p_2$, increasing $|p_2|$ by less than $b$, so the adjusted side is a full power $B^n$. Move each cut endpoint by less than $b$ to a $B$-phase vertex. Let the resulting cut labels be $T,S$, each of length less than $(\gamma^{-1}+1)b$. Let $P$ be the unchanged label of $p_1$ and $Q$ the adjusted label of $p_2$. If $B^{r_1}$ is the initial part ending at the $T$-endpoint and $B^{r_2}$ the terminal part beginning at the $S$-endpoint, the cut-disc equations are
\[
 VQB^{r_1}T^{-1}=1,\qquad PWSB^{r_2}=1.
\]
Using $V=W$, obtain
\begin{equation}\label{longshortidentity}
 (PT)B^{-r_1}(Q^{-1}S)=B^{-r_2}.
\end{equation}
The two connector factors have total length less than $K_0b$. Length comparison on the cut discs gives, for $j=1,2$,
\begin{align*}
 r_j&>c_0L_0/b-K_1\\
 &>\frac{c_0\gamma}{2+\gamma}(c_0N_*/2-2\alpha)-K_1
 \ge\gamma N_*/8.
\end{align*}
Here $L_0>\gamma|q_1|/(2+\gamma)$ and $|q_1|>(c_0N_*/2-2\alpha)b$, and the final inequality is \eqref{Q2}. In particular both powers are nonempty. By \eqref{Q3}, Lemma~\ref{pcm}(b), with root $B^{-1}$, applies to \eqref{longshortidentity}. Thus $PT=B^e$.

The remaining prefix-disc equation is $PATB^{n-r_1}=1$, so
\[
 A=P^{-1}B^{r_1-n-e}P.
\]
Simple-power conjugacy makes the exponent $\pm1$, proving both assertions.
\end{proof}

\subsection{The comparison required at the next rank}
\begin{lemma}\label{weakcomparison}
Let $A,B$ be current-simple. Suppose a reduced disc has contour $p_1q_1p_2q_2$ with $q_1,q_2$ periodic in $A,B$, respectively, and
\[
 |q_2|\ge\eps n_0|B|,\qquad
 |p_1|,|p_2|<\zeta n_0\min(|A|,|B|).
\]
Then either $|q_1|<(1+\gamma)|A|$, or $A\sim B^{\pm1}$. If $A$ and $B$ coincide literally up to inversion, they must have opposite orientations and the two sections have compatible phase vertices.
\end{lemma}
\begin{proof}
Put $a=|A|$, $b=|B|$, and $s=\min(a,b)$. Trim with loss at most $T_0s$. Length comparison, \eqref{Q4}, and $\zeta\ll\eps$ ensure that both sides satisfy the trimming hypothesis. The retained $B$-side is longer than $N_*b/2$, and the new connectors are shorter than $\alpha s$.

If the retained $A$-side is longer than $(1+\gamma/2)a$, apply Lemma~\ref{longshort}. If $A=B^{\pm1}$ literally, the retained sides are both longer than $h_ga$ by \eqref{Q1} and length comparison. Lemma~\ref{orientation} excludes equal orientations, and Lemma~\ref{pcm}(a) gives compatibility for opposite orientations.

In the other case, comparison with the retained $B$-side gives
\[
 (c_0N_*/2-2\alpha)b<(1+\gamma/2)a.
\]
The discarded $A$-length is at most
\[
 T_0s\le T_0b
 <\frac{T_0(1+\gamma/2)}{c_0N_*/2-2\alpha}\,a<\gamma a/2
\]
by \eqref{Q5}. Restoring it gives $|q_1|<(1+\gamma)a$. When $A=B^{\pm1}$, this short alternative is itself incompatible with the original $B$-side length and (G2), so the orientation and compatibility conclusion always follows in that case.
\end{proof}

\subsection{Completion of the rank induction}
The preceding arguments use arbitrary positive cells only through Condition A and smoothness. The only removal that uses a power label is the explicitly described surgery on a selected-period boundary. In particular no assumption that the inserted cells are powers has been made.

Select the next periods. They are simple at every earlier rank: failure of simplicity in an earlier quotient persists in all later quotients before their selection. If two selected periods are conjugate up to inversion below both defining ranks, unequal lengths contradict simplicity and equal lengths contradict the choice of representatives. This is the period-recognition assertion used in Section~\ref{insertionsection}.

For a power--power strip in the next diagram induction, rank decrease puts the strip below both defining ranks. Lemma~\ref{weakcomparison} applies there. Its long alternative identifies the roots, and period recognition and its phase conclusion give precisely a compatible pair. A self-contiguity is excluded by its orientation clause. Thus only a strictly earlier stage is used in the next $D_N,S_N$ proof.

The dependence order is therefore: earlier-rank comparison; current $D_N,S_N$ for all $N$; root and torsion facts; Lemma~\ref{pcm} in increasing root-length sums; Lemmas~\ref{orientation}--\ref{weakcomparison}; next period selection. Starting with the old ranks, it proves these assertions at $\rho$ and then at every subsequent integer rank. This completes the proof of Proposition~\ref{insertion}.

\section{Proof of Theorem~\ref{main}}
\begin{proof}
 Assume that the chosen tuple satisfies all conclusions of
Lemma~\ref{goodevent}. Let $H_\infty$ be the group given by the reselected presentation. Every defining relator of $H_\infty$ holds in $Q_\ell$, because it is either an $M$th power or one of the $R_i$.

Conversely, all words of length at most $d$ have $M$th power one in $H_d=G_d$. Inductively, consider a word $W$ of length $j>d$ just before rank $j$. If it is simple, it is conjugate up to inversion to one of the selected rank-$j$ representatives, whose $M$th power is imposed. Otherwise it is trivial or conjugate to a power of a shorter word or a previously selected period; those elements already have $M$th power one. Thus every word of length at most $j$ has $M$th power one in $H_j$. It follows that $H_\infty$ has exponent dividing $M$, and it also kills every $R_i$. The two maps fixing the original generators therefore identify
\[
 H_\infty=B/\normal{R_1,\ldots,R_k}=Q_\ell.
\]
The last equality uses conjugacy in $B$, and does not require a chosen conjugator in a modified stage.

If a word $w$ represents a nontrivial element of $B$ killed in $Q_\ell$, any reduced diagram for $w$ must contain a cell of rank at least $\rho$. Otherwise the diagram is old and $w=1$ in $B$. Proposition~\ref{insertion} and the perimeter lemma give
\[
 |w|>(1-\beta)n_0\rho.
\]
By \eqref{insertionscale}, this is larger than $\kappa\ell$ for all sufficiently large $\ell$, for example with $\kappa=(1-\beta)b/4$. The probability estimate follows from Lemma~\ref{goodevent}.

By the restricted Burnside theorem there is a finite bound $F(m,M)$ on the order of any finite $m$-generated group of exponent dividing $M$ \cite{Zelmanov}. Since $B$ is infinite, choose a fixed $R$ with $|B_B(R)|>F(m,M)$. Once $2R<\kappa\ell$, the kernel estimate makes the projection injective on $B_B(R)$. Thus $Q_\ell$ has more than $F(m,M)$ elements and cannot be finite.
\end{proof}
\section{The Gromov density model}\label{gromovsection}
Put $q=2m-1$, so that $H=\log q$, and let $\mathcal C_\ell$ be the set of cyclically reduced words of length $\ell$ in the free generators and their inverses. For a fixed density $d\in[0,1)$, put
\[
 K_\ell=\lfloor q^{d\ell}\rfloor.
\]
Choose $W_1,\ldots,W_{K_\ell}$ independently and uniformly from $\mathcal C_\ell$. This is the version of Gromov's density model in which sampling is with replacement; the version using a uniform $K_\ell$-element subset has the same conclusion below. See \cite{Ollivier} for the density model. Write
\[
 G_{\ell,d}=F_m/\normal{W_1,\ldots,W_{K_\ell}},
 \qquad
 Q_{\ell,d}^{(M)}
 =B_m(M)/\normal{\overline W_1,\ldots,\overline W_{K_\ell}},
\]
where $\overline W_i$ is the image of $W_i$ in $B_m(M)$. Equivalently, $Q_{\ell,d}^{(M)}$ is the quotient of $G_{\ell,d}$ by all $M$th powers.

\begin{theorem}[Positive density in the Burnside variety]\label{gromovtheorem}
For every $m\ge2$ there are $d_0>0$, $M_0$ and $\kappa>0$ with the following property. For every fixed odd $M\ge M_0$ and every fixed $0\le d< d_0$, there is $c>0$ such that, for all sufficiently large $\ell$,
\[
 \Pr\!\left\{
 \ker\bigl(B_m(M)\longrightarrow Q_{\ell, d}^{(M)}\bigr)
       \cap B_{B_m(M)}(\kappa\ell)=\{1\}
 \right\}\ge1-e^{-c\ell}.
\]
In particular, $Q_{\ell, d}^{(M)}$ is infinite with probability at least $1-e^{-c\ell}$. The conclusion also holds for independent uniform freely reduced words of length $\ell$, and for a uniformly chosen $K_\ell$-element subset of $\mathcal C_\ell$.
\end{theorem}

Fix the parameters of Section~\ref{parametersection}, depending only on $m$, and put
\begin{equation}\label{densitysaving}
 \sigma=\frac14\min\{(1-b)(h_*-D/2),\,h_*\tau\}>0.
\end{equation}
One admissible density bound is
\begin{equation}\label{densitythreshold}
 \delta_0=\frac{\sigma}{4H}
 =\frac{\min\{(1-b)(h_*-D/2),\,h_*\tau\}}{16H}.
\end{equation}
This is a sufficient positive bound, not an assertion about the critical density. We enlarge $M_0$ below to impose one additional cogrowth estimate, keeping all the preceding parameters fixed.

\subsection{Exceptional tuples in a  ball in $B_m(M)$}
Let $\mathcal E_1(\ell)\subset B_B(\ell)$ consist of the elements whose minimum conjugacy length is at most $b\ell$, together with the elements violating one of the one-core exclusions (i), (iii), or (iv) of Lemma~\ref{goodevent}. Let $\mathcal E_2(\ell)\subset B_B(\ell)^2$ consist of the ordered pairs violating its two-core exclusion (ii). A violation means that at least one choice of shortest representatives, signs, positions, and the specified auxiliary data is a forbidden witness. Direct natural self-identifications remain allowed, exactly as in Lemma~\ref{goodevent}.

\begin{lemma}\label{ballexceptions}
For each fixed admissible $M$, with $h=h_M$,
\begin{equation}\label{ballexceptionbound}
 |\mathcal E_r(\ell)|\le
 \exp\bigl(rh\ell-2\sigma\ell+o(\ell)\bigr),
 \qquad r\in\{1,2\}.
\end{equation}
The bounds concern one element or one ordered pair and are independent of the number of sampled relators.
\end{lemma}
\begin{proof}
For a nontrivial conjugacy class with minimum length $L\le\ell$, sum \eqref{fiber} over the radii from $L$ to $\ell$. This gives
\begin{equation}\label{ballfiber}
 |\operatorname{Conj}(U)\cap B_B(\ell)|
 \le\operatorname{poly}(\ell)
 \exp\bigl(D(\ell-L)/2+9\gamma_cD\ell\bigr).
\end{equation}
The identity class is a singleton and also satisfies this upper bound when $L=0$. Using $s(L)\le\exp(hL+E(\ell))$, the number of elements with minimum conjugacy length at most $b\ell$ is therefore at most
\[
 \exp\bigl(h\ell-
       [(1-b)(h-D/2)-9\gamma_cD]\ell+o(\ell)\bigr).
\]
By \eqref{codingchoice} and $h\ge h_*$, the saving in square brackets is greater than $\tfrac12(1-b)(h_*-D/2)\ge2\sigma$.

For any of the compressed descriptions in (i)--(iii), with $r$ core lengths $L_1,\ldots,L_r$, the code count is at most
\[
 \exp\bigl(h(\textstyle\sum_i L_i-\tau\ell)+o(\ell)\bigr).
\]
Multiply by \eqref{ballfiber} for each decoded class. Since
\[
 h\sum_iL_i+\frac D2\sum_i(\ell-L_i)\le rh\ell,
\]
the resulting count is at most
\[
 \exp\bigl(rh\ell-[h_*\tau-9r\gamma_cD]\ell+o(\ell)\bigr).
\]
The same calculation for the literal-root scheme (iv) uses \eqref{literalbad} to give the saving $h_*\tau-9\gamma_cD$. In all these cases the saving is greater than $\tfrac34h_*\tau\ge3\sigma$, by \eqref{codingchoice}. The finitely many schemes, length vectors, phases and other polynomial choices contribute $o(\ell)$; the entropy errors are bounded multiples of $E(2\ell)=o(\ell)$. This proves \eqref{ballexceptionbound}. As before, the codes count the existence of a witness without enumerating all geodesic spellings.
\end{proof}

\subsection{Transfer from group elements to random words}
The required additional input is the following consequence of Adian's cogrowth estimate \cite{AdianRW}; its adjustable parameter is explained in \cite[Remark 8.4]{GGMT}. For every fixed $m\ge2$ and every $\eta>1/2$, all sufficiently large odd $M$ satisfy
\begin{equation}\label{cogrowthinput}
 T_n:=\#\{W:\ W\text{ is freely reduced},\ |W|=n,
                     \ \overline W=1\text{ in }B_m(M)\}
 \le C q^{\eta n}\qquad(n\ge0),
\end{equation}
where $T_0=1$ and $C$ may depend on $m,M,\eta$. More precisely, the upper bounds for the cogrowth rate tend to $q^{1/2+2/Q}$ as the odd exponent tends to infinity, and the auxiliary parameter $Q$ can be chosen arbitrarily large. Choose $Q$ with $1/2+2/Q<\eta$, and then enlarge the exponent threshold. The strict inequality for the limiting rate yields \eqref{cogrowthinput}, after increasing $C$ to cover finitely many lengths. Thus this input holds for every sufficiently large odd exponent, fixed before $n$ varies.

\begin{lemma}\label{spellingmeasure}
Suppose \eqref{cogrowthinput} holds with $1/2<\eta<1$. Let $\mu_\ell$ be the distribution on $B$ obtained by evaluating a uniform freely reduced word of length $\ell$. Then
\begin{equation}\label{collisionbound}
 \sum_{g\in B}\mu_\ell(g)^2
 \le C' q^{-2(1-\eta)\ell}.
\end{equation}
The same bound, with another constant, holds for the distribution $\nu_\ell$ obtained from a uniform cyclically reduced word of length $\ell$. Consequently, for $r\in\{1,2\}$ and $E\subset B_B(\ell)^r$,
\begin{equation}\label{spellingtransfer}
 \nu_\ell^{\otimes r}(E),\ \mu_\ell^{\otimes r}(E)
 \le C_r |E|^{1/2}q^{-r(1-\eta)\ell}.
\end{equation}
\end{lemma}
\begin{proof}
Count pairs of freely reduced words $U,V$ of length $\ell$ with the same value in $B$. If their maximal common suffix is $S$, of length $j$, write $U=PS$ and $V=QS$. The word $PQ^{-1}$ is freely reduced, has length $2(\ell-j)$, and represents the identity. Conversely, $PQ^{-1}$ determines $P$ and $Q$ by splitting it into equal halves. There are at most $C_mq^j$ possible suffixes, including the empty suffix. Hence the number of pairs is at most
\[
 C_m\sum_{j=0}^{\ell}q^jT_{2(\ell-j)}
 \le C''q^{2\eta\ell}
       \sum_{j=0}^{\ell}q^{(1-2\eta)j}
 \le C'''q^{2\eta\ell}.
\]
This includes $U=V$, for which $j=\ell$. Dividing by the square of the number $2m q^{\ell-1}$ of freely reduced words proves \eqref{collisionbound}.

For $\ell\ge2$, at most a fraction $1/q$ of freely reduced words fail to be cyclically reduced: after the first $\ell-1$ letters there is at most one last letter causing this failure. Thus the conditioning event has probability at least $p_m=1-1/q>0$, and
\[
 \nu_\ell(g)\le p_m^{-1}\mu_\ell(g)
 \qquad(g\in B).
\]
This proves the cyclically reduced version of \eqref{collisionbound}. Finally apply Cauchy--Schwarz to the sum over $E$, using
\[
 \sum_{(g_1,\ldots,g_r)\in B^r}
       \prod_{i=1}^r\mu_\ell(g_i)^2
 =\left(\sum_{g\in B}\mu_\ell(g)^2\right)^r,
\]
and the analogous identity for $\nu_\ell$.
\end{proof}

\subsection{The density estimate and infiniteness}
\begin{proof}[Proof of Theorem~\ref{gromovtheorem}]
Choose
\[
 \eta=\frac12+\frac{\sigma}{8H}
\]
and enlarge $M_0$ so that \eqref{cogrowthinput} holds for all odd $M\ge M_0$. This choice satisfies $1/2<\eta<1$. Combining Lemmas~\ref{ballexceptions} and \ref{spellingmeasure}, and using $h\le H$, gives, for $r=1,2$,
\begin{align*}
 \nu_\ell^{\otimes r}(\mathcal E_r(\ell))
 &\le\exp\left(
       \left[\frac{rh}{2}-\sigma-r(1-\eta)H\right]\ell
        +o(\ell)\right)\\
 &\le\exp\bigl(-(1-r/8)\sigma\ell+o(\ell)\bigr)\\
 &\le\exp\bigl(-3\sigma\ell/4+o(\ell)\bigr).
\end{align*}
The same estimates hold with $\mu_\ell$ in place of $\nu_\ell$. Only one or two sampled words are involved in each test; in particular, we have not conditioned the entire exponentially long tuple on cyclic reduction.

There are $K_\ell$ one-core tests and at most $K_\ell^2$ ordered two-core tests. Let $\mathsf{Bad}_\ell$ be the event that at least one of these tests fails. A union bound gives
\begin{equation}\label{densityfailure}
 \Pr(\mathsf{Bad}_\ell)
 \le\exp\bigl((2d H-3\sigma/4)\ell+o(\ell)\bigr)
 \le\exp\bigl(-\sigma\ell/4+o(\ell)\bigr),
\end{equation}
where the last inequality uses \eqref{densitythreshold}. The constants and the $o(\ell)$ terms in the one-core and two-core estimates do not depend on the indices of the sampled words.

On the complementary event, choose any cyclic core $R_i$ of each $\overline W_i$. Then $b\ell<|R_i|\le\ell$, and all the deterministic exclusions in Lemma~\ref{goodevent} hold. Proposition~\ref{insertion} applies with $k=K_\ell$. To make the uniformity here explicit, its diagram induction is on rank and on the number of cells in an individual diagram, not on the size of the list of defining relators. Every mixed-strip test involves one or two members of that list and is covered by the exclusions just established, including the self-contact tests. Its geometric inequalities, rank decrease and compatible-pair replacements contain no factor depending on $k$. At each fixed $\ell$ the list is finite, and the graded type still has one coordinate for each rank. Thus the same induction proves the proposition for this list, with the same parameters. The fixed-$k$ assumption was used only in the probability union bound, which has now been replaced by \eqref{densityfailure}.

The identification of the reselected presentation with the quotient, and the perimeter argument in the proof of Theorem~\ref{main}, now give
\[
 \ker(B\longrightarrow Q_{\ell, d}^{(M)})
     \cap B_B(\kappa\ell)=\{1\},
 \qquad \kappa=(1-\beta)b/4.
\]
These steps only use the lengths and exclusions for the chosen cores; they do not require $|\overline W_i|_B=\ell$. For infiniteness, the restricted Burnside theorem gives a bound $F(m,M)$ on the order of a finite $m$-generated group of exponent dividing $M$. Choose a fixed $R$ with $|B_B(R)|>F(m,M)$. Once $2R<\kappa\ell$, the kernel estimate makes the quotient map injective on $B_B(R)$, so the quotient cannot be finite. This proves the assertions for independent sampling.

Finally, \eqref{tauchoice} and \eqref{densitythreshold} imply $d_0<1/2$. For independent uniform samples from $\mathcal C_\ell$, the probability that two sampled words coincide is at most
\[
 \binom{K_\ell}{2}|\mathcal C_\ell|^{-1}
 \le C_m q^{(2 d-1)\ell}=o(1).
\]
Conditional on no repetitions, their unordered set is a uniform $K_\ell$-element subset. Dividing the exponentially small failure probability by the probability of no repetitions preserves exponential decay. This proves the final version of the theorem.
\end{proof}
\subsection{The density treshhold}
Write $T_n$ for the number of freely reduced words of length $n$ representing
$1$ in $B$, and define
\[
 \eta_M=\lim_{\substack{n\to\infty\\ n\ \mathrm{even}}}
 \frac{\log_q T_n}{n}.
\]
The standard cogrowth facts give existence of this limit and
$\eta_M>1/2$, since $B$ is a proper quotient of $F_m$
\cite[Section 1.2]{Ollivier}. Adian's estimates imply
\begin{equation}\label{cogrowthbounds}
 \frac12<\eta_M<1,
 \qquad
 \lim_{\substack{M\to\infty\\ M\ \mathrm{odd}}}\eta_M=\frac12.
\end{equation}
For the second assertion, one uses the adjustable auxiliary parameter in
Adian's estimate \cite{Adian}; see \cite[Remark 8.4]{GGMT}. In particular,
along even $n$,
\begin{equation}\label{count}
 T_n=q^{\eta_M n+o(n)}.
\end{equation}

The following is a specialization of Ollivier's high-density argument
\cite[Section 5.2]{Ollivier}, which does not require hyperbolicity of the base
group. The direct proof makes the reduction and parity conditions explicit.

\begin{proposition}\label{collapse}
If $d>1-\eta_M$, there is $c>0$ such that, for all sufficiently large $\ell$,
\[
 \Prob\{Q_{\ell,d}^{(M)}\ne1\}\le m\exp(-\exp(c\ell)).
\]
The same conclusion holds for uniform freely reduced relators and for
sampling a uniform $K_\ell$-element subset of either set of words.
\end{proposition}

\begin{proof}
Fix $a\in X$ and choose $b\in X\setminus\{a\}$. For each sufficiently large
$\ell$, put
\[
 t=\begin{cases}1,&\ell\text{ odd},\\ M+1,&\ell\text{ even},\end{cases}
 \qquad n=\ell-4-t.
\]
Then $n$ is positive and even. Let $u$ be any freely reduced identity word
of length $n$ in $B$, with first letter $f$ and last letter $z$. Choose a
letter $s\in X^{\pm1}$ avoiding $b^{-1},f^{-1},z$, using a fixed ordering to
make the choice deterministic. Such a letter exists because $2m\ge4$.
Set $v=bs$ and
\[
 w(u)=vuv^{-1}a^t=bsu s^{-1}b^{-1}a^t.
\]
There is no free cancellation at any junction. The first letter is $b$
and the last is $a$, so $w(u)$ is also cyclically reduced. Its length is
$\ell$, and in $B$ it represents $a$, because $u=1$ and $t\equiv1\pmod M$.
The map $u\mapsto w(u)$ is injective: $u$ occupies the fixed positions
$3,\ldots,n+2$. The added length is at most $M+5$, a constant because $M$
is fixed.

Consequently, if $p_{a,\ell}$ is the probability that a uniform word in
$\mathcal C_\ell$ represents $a$, then \eqref{count} and
$|\mathcal C_\ell|\le2m q^{\ell-1}$ give
\[
 p_{a,\ell}\ge\frac{T_n}{2m q^{\ell-1}}
       =q^{-(1-\eta_M)\ell+o(\ell)}.
\]
The same lower bound holds when sampling all freely reduced words.
Thus the probability that none of the $K_\ell$ relators represents $a$ is
at most
\[
 (1-p_{a,\ell})^{K_\ell}
 \le\exp(-K_\ell p_{a,\ell})
 \le\exp\!\left(-q^{(d+\eta_M-1)\ell+o(\ell)}\right).
\]
Since $d+\eta_M-1>0$, this is at most $\exp(-\exp(c\ell))$ for some
$c>0$, uniformly over the finitely many $a\in X$. If every generator is
represented by a sampled relator, every generator is killed in the quotient.
A union bound proves the assertion.

For sampling without replacement, let $N$ be the number of admissible words
and $A$ the number representing $a$. The probability of missing all of them is
\[
 \frac{\binom{N-A}{K_\ell}}{\binom{N}{K_\ell}}
 \le(1-A/N)^{K_\ell}\le\exp(-K_\ell A/N),
\]
with numerator zero when $K_\ell>N-A$. The same estimate applies.
\end{proof}

\subsection{What remains to be proved}
Theorem 8.1 gives, for each $m$,
constants $d_0>0$ and $M_0$ such that for every fixed odd $M\ge M_0$ and every
fixed $d<d_0$, the group $Q_{\ell,d}^{(M)}$ is infinite with probability at
least $1-e^{-c'\ell}$. It also gives a linear lower bound for the length of a
nontrivial kernel element of $B\to Q_{\ell,d}^{(M)}$. Taking that theorem as
an input, and decreasing $d_0$ if necessary, we obtain the following picture:
\begin{center}
\begin{tabular}{ll}
\toprule
Density & Conclusion as $\ell\to\infty$\\
\midrule
$0\le d<d_0$ & Infinite with probability tending to one\\
$d_0\le d\le1-\eta_M$ & Not settled by these arguments\\
$d>1-\eta_M$ & Trivial with probability tending to one\\
\bottomrule
\end{tabular}
\end{center}
These bounds do not establish a sharp transition. Monotonicity under adding
relators alone does not close the gap.

\begin{question}
For each fixed sufficiently large odd $M$, is $Q_{\ell,d}^{(M)}$ infinite
asymptotically almost surely for every $d<1-\eta_M$?
\end{question}

A positive answer would identify the critical density as
$d_c(m,M)=1-\eta_M$; no assertion at equality is included. This is a candidate
threshold, not a consequence of the present estimates. By
\eqref{cogrowthbounds}, the candidate tends to $1/2$ as odd $M\to\infty$.
It is strictly smaller than $1/2$ for every fixed $M$.
Thus infiniteness for all $d<1/2$ with one fixed exponent is impossible.
The different statement
\[
 \forall d\in[0,\tfrac12)\ \exists M_0(m,d)\ \forall\text{ fixed odd }M\ge M_0(m,d),
 \qquad \Prob\{Q_{\ell,d}^{(M)}\text{ infinite}\}\longrightarrow1
\]
is compatible with the collapse bound and would follow from a positive
answer to the question.
\section{Acknowledgement} The paper was written with Open AI Chat GPT 6 Pro editing assistance.  It also helped with selection of parameters and computation of estimates. All was carefully checked by the
authors and the authors carry the full responsibility for the content of the paper.

\end{document}